\documentclass[a4paper]{article}
\usepackage{times,cite,w-thm}
\usepackage[T1]{fontenc}
\usepackage[utf8]{inputenc}
\theoremstyle{plain}

\usepackage{graphicx}
\usepackage{amsmath}% for pmb

\newcommand{\bx}{{\bf{x}}}
\newcommand{\bu}{{\bf{u}}}
\newcommand{\bv}{{\bf{v}}}

\newcommand{\bsigma}{{\boldsymbol{\sigma}}}

\newcommand{\tr}{{\rm tr}}

\newcommand{\dev}{{\rm dev}}

\renewcommand{\div}{{\rm div}}

\newcommand{\N}{{\rm I\kern-.25em N}}
\newcommand{\D}{{\rm I\kern-.25em D}}

\newcommand{\bepsilon}{\mbox{\boldmath$\varepsilon$}}
\newcommand{\btau}{\mbox{\boldmath$\tau$}}

\newcommand{\bF}{{\bf F}}

\newcommand{\bp}{{\bf p}}
\newcommand{\bq}{{\bf q}}

\newcommand{\bm}{{\bf m}}
\newcommand{\be}{{\bf e}}
\newcommand{\bI}{{\bf I}}

\newcommand{\bW}{{\bf W}}

\newcommand{\bR}{{\bf R}}
\newcommand{\cF}{{\cal F}}
\newcommand{\cG}{{\cal G}}

\begin{document}
%% \def\leftmark{Session  }
%%
%%    The information for the title page will be placed between
%%    \begin{document} and \maketitle. The order of most entries
%%    is determined by the class file and can not be changed by
%%    rearranging them. The maketitle command follows after the
%%    abstract.
%%
%%    Most of the following commands will be completed by the publisher.
%%
%%    \renewcommand{\copyrightyear}{2016}
%%    \DOIsuffix{pamm.20161zzzz}
%%    \Volume{16} 
%%    \Year{2016} 
%%    \pagespan{1}{}
%%
%%    The short title is optional:

%\TitleLanguage[EN]
\title{Least Squares Finite Element Methods for Sea Ice Dynamics}
\date{5.09.2018}
\author{Fleurianne Bertrand}
%% Please do not enter footnotes or \inst{}-notes into the optional
%% argument of the author command. 

%% Please delete not needed author entries.
%% Information for the first author.
%\author{\firstname{Fleurianne}  \lastname{Bertrand}\inst{1,}%
%  \footnote{Corresponding author: e-mail \ElectronicMail{fleurianne.bertrand@uni-due.de}}}

%\address[\inst{1}]{\CountryCode[DE]Thea Leymann Stra\ss e 9, 45127 Essen }

%\AbstractLanguage[EN]

%% maketitle must follow the abstract.
\maketitle                   % Produces the title.

\begin{abstract}
A first-order system least squares formulation for the sea-ice dynamics is presented.  In addition to the displacement field, the stress
tensor is used as a variable. As finite element spaces, standard conforming piecewise polynomials for the
displacement approximation are combined with Raviart-Thomas elements for the rows in the stress
tensor.  Computational results for a test problem illustrate the least-squares approach.
\end{abstract}

\section{Introduction}
Ice and snow covered surfaces reflect more than half of the solar radiation they are recieving and play therefore a major role in climate modelling. Each year, Antarctic sea ice extent reaches its maximum (17-20 million square kilometers) in September and 
 its minimum (3-4 million square kilometers) in February. These important oscillations make the current predictive models of Antarctic sea ice require an accurate knowledge and understanding of the processes. Developing computational sea-ice modelling based on observed and measured data to study and predict the break-up and fracture evolution of sea-ice during the Antarctic spring was one of main scientific aims of the 
% ToDo CITE NASA
 Winter 2017 cruise (Voyage 25) of the S.A. Agulhas II. This was funded by DST/NRF and took place
from 28 June to 13 July 2017. %The cruise consisted of 3 legs with several national and international
%participants, including the author. %Details on the framework in which this voyage
%originated are illustrated in the final Chapter 15.
%The major goal of the cruise was to capture the winter recharge conditions in the Atlantic-Indian
%Ocean sector of the Southern Ocean, from open ocean conditions to the marginal ice zone (MIZ), with
%a special emphasis on documenting the variability of the MIZ. The voyage was also endorsed by
%GEOTRACES as a process study and it covered the sampling of the WOCE I06 transect. The cruise
%participants are listed in the table below and portrayed in Figure 1 on their way back to Cape Town.
%There were 11 scientific teams on board (Chapters 3-13), with the addition of an artist that made an
%alternative documentation record of the cruise that brought art and science together (Chapter 14).
%The were multiple 
%To document sea ice growth processes, physical properties (temperature, porosity, salinity,
%solid-fluid-gas volume fractions), and the mechanical properties (strength, stiffness, fracture
%toughness, viscous-elasticity) of sea ice (Resp: MacHutchon, UCT)
%
%To develop computational sea-ice modelling based on observed and measured data to study and
%predict the break-up and fracture evolution of sea-ice during the Antarctic spring; advance
%understanding of the complex material behaviour of sea-ice in terms of its liquid and gas-filled
%porous structural composition; elucidate the mechanisms of biogeochemical exchanges with the
%underlying ocean and the brine expulsion (Resp: Skatulla UCT, Ricken UniD-E, Bertrand TUD)

Sea ice is a complex material which is formed by the freezing of sea water. %During its formation and
%growth, the sea ice structure is profoundly modified by the interaction of physical, biological and
%chemical processes, and becomes a heterogeneous semi-solid matrix in its simplest form (Hunke et al.,
%2010). These ice processes affect the climate (Light, 2003), ocean mixed layer structure (Morison and
%Smith, 1981), biological activity in ice and water (Cota, 1985), and the ocean heat flux at the underside
%of the ice (Maykut and McPhee, 1995). In winter, the Southern Ocean sea ice forms a barrier mediating
%the exchange of heat and momentum between the atmosphere and the ocean. The gradual ice break-up
%in spring and summer renders the barrier permeable and influences the ice dynamics in the region
%(Langhorne et al., 1998).
Since the ice stress is a source in the other equations of the climate models, its approximation plays an important role in the simulations of the ice. They can be computed from the velocity in a post-processing step, but the loss of accuracy due to the reconstruction step can lead
to non-physical solutions. An alternative approach consists in the use of variational formulations involving
the stress $\bsigma\in H(\div,\Omega)$ as an independent variable. Appropriate finite element spaces based on a triangulation $\mathcal T$ are the $H(\div,\Omega) $-conforming spaces, e.g. the Raviart-Thomas Space.

\section{Problem Formulation}
As most sea ice dynamic models currently used, our model is based on the viscous-plastic formulation introduced by Hibler \cite{Hibler:79}. There, sea ice is modeled by its velocity $\bu$, the ice concentration $A$ and the average ice height $H$ over a domain $ \Omega$. % TODO {\color{red}Constraints???}
The model consists in a momentum equation for the velocity $\bu$ and the balance laws  for ice concentration $A$ and the average ice height $H$. Neglegting the thermodynamical effects, i.e. the source terms in these balance laws, the model can be written as
\begin{align}\label{strong}\begin{split}
\rho_{ice} H \frac{\partial \bu}{\partial t} + \bF(\bu) - \div\ {\bsigma(\bu,A,H)}&=0, \\
\frac{\partial A}{\partial t}+\div(\bu A)&=0 , \quad
\frac{\partial H}{\partial t}+\div(\bu H)=0 \ ,
\end{split}\end{align}
where the force term involving the ice, air and water densities $\rho_{ice}$,$\rho_a$ and $\rho_o$, the air and water drag coefficients $C_a$ and $C_o$, the coriolis parameter $f_c$, the radial unit vector $\be_r$ and the velocity fields $ \bv_o$ and  $\bv_a$ of ocean and atmospheric flow is given by 
% TODO Epsilon u
\begin{align}
\bF(\bv) = f_c\be_r \times (\bv-\bv_o) - 
\underbrace{\rho_a C_a \|\bv_a\|_2\bv_a}_{=:\btau_a} - 
\underbrace{\rho_o C_o \|\bv_o-\bv\|_2(\bv_o-\bv)}_{=:\btau_o(\bv)}
\end{align}
and the stress-strain relation involving the ice strength parameter $P^\star$ and the ice concentration parameter $C$ is given by 
\begin{align}\label{stressstrain}\begin{split} 
&\bsigma = \frac P 2 
	\left(
		\frac { {\dev \ \bepsilon(\bu)}+  2 {\tr\ \bepsilon(\bu)}  \bI} {\Delta(\bu)}  -\bI 
	 \right) \quad
 \text{with } P = P^\star H e^{-C(1-A)}\\
& \text{and }\ \Delta(\bu) = \sqrt{\dev\ \bepsilon(\bu):\dev\ \bepsilon(\bu)+{4}\tr(\bepsilon(\bu))^2+\Delta_{min} ^2}
% -\frac{P}{2\sqrt{\frac{1}{2}(\dev\ \bepsilon(\bu)):(\dev\ \bepsilon(\bu))+\frac{3}{4}\tr(\bepsilon(\bu))}}
\end{split}\end{align}
where $\Delta_{min} = 2 \cdot 10^{-9}\ s^{-1}$ is a limitation for $\Delta(\bu)$.
In \cite{MehRic:17}, the authors propose a variational formulation where $(\bu, \bp)$ with $\bp = (A,H)$ is sought in $\left(H^1_0(\Omega)\right)^2 \times \left(L^2(\Omega)\right)^2$ such that
\begin{align}\begin{split}
\left( \rho_{ice} H \frac{\partial \bu}{\partial t},\bv\right)+(\bF(\bu),\bv)+(\bsigma(\bu,H,A),\nabla \bv) & = 0, \\
\left(\frac{\partial \bp}{\partial t}+\nabla \bp \cdot \bu+\div(\bu)\bp,\bq\right)  & = 0
%\left(\frac{\partial A}{\partial t}+\bu\cdot\nabla A+A\div(\bu),p\right)  &= 0\\ 
%\left(\frac{\partial H}{\partial t}+\bu\cdot\nabla H+H\div(\bu),q\right) &=0 
\end{split}\end{align}
holds for all $(\bv,\bq) \in \left(H^1_0(\Omega)\right)^2 \times \left(L^2(\Omega)\right)^2$.
The constraints $H \geq 0$ and $A \in [0,1]$ are embedded in the trial-spaces and are realized by a projection of the solution. %The purpose of the next section is to present a stress-based mixed finite element method {\color{red}TODO} variational inequalities.
\begin{figure}[t]
\includegraphics[width=0.5\textwidth]{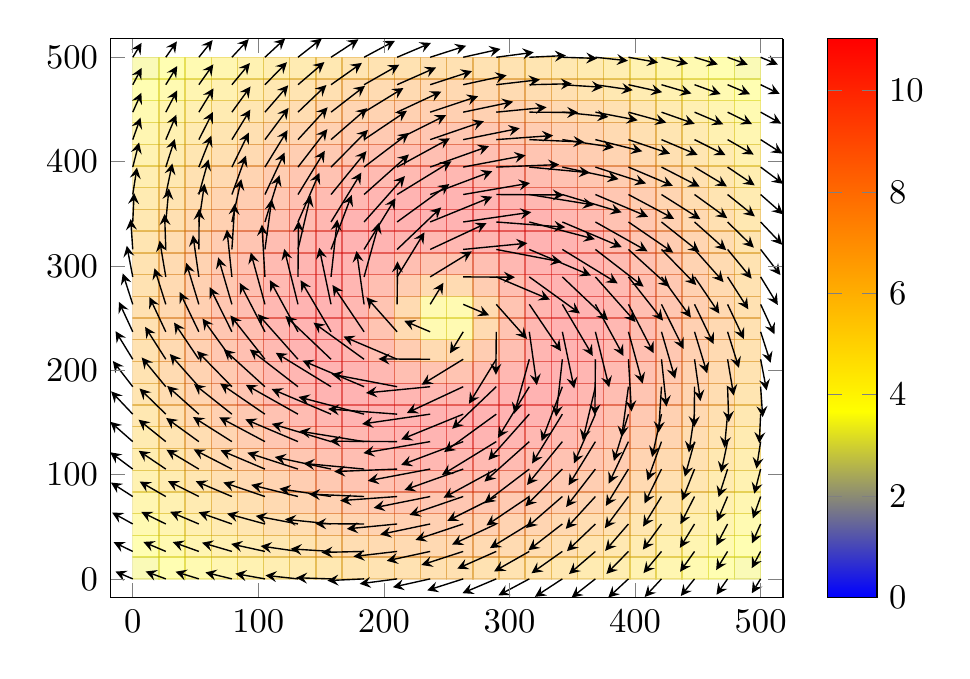}
\includegraphics[width=0.5\textwidth]{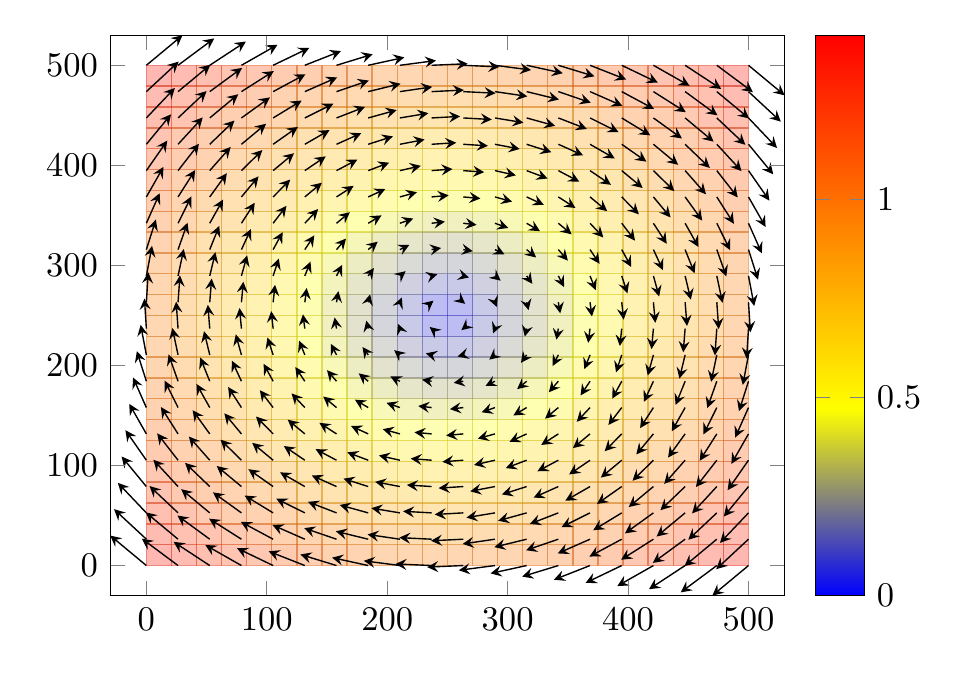}
\caption{Wind field at $t=0$ (left) and Ocean current (right)}
\label{fig:3}
\end{figure}

\section{A Least-Squares Method}
The Least-Squares Mehtod (see \cite{BocGun:09}) consists in minimizing the $L^2$-residuals in the partial differential equations. Therefore, we insert define a new variable $\bsigma$ for the stress and consider the stress-strain relationship \eqref{stressstrain}  as and additional equation in order to obtain the following first order system for $(\bsigma,\bu,A,H)$:
\begin{align*}\label{fos}
\rho_{ice} H \frac{\partial \bu}{\partial t} + \bF(\bu) - \div\ {\bsigma}&=0  & \frac{\partial A}{\partial t}+\div(\bv A)&=0 \\
 \frac {P(A,H)} 2 \left( \frac {\dev \ \bepsilon(\bu)} {\Delta(\bu)} + \frac  {2\tr\ \bepsilon(\bu)} {\Delta(\bu)} \bI -\bI \right) &=\bsigma & \frac{\partial H}{\partial t}+\div(\bv H)&=0 \
\end{align*}
The least-squares functionals then reads \begin{align}
\cF (\bsigma,\bu,H) = \cF_m (\bsigma,\bu,A,H)+\cF_c (\bsigma,\bu,A,H)+\cF_e (\bsigma,\bu,A,H)\end{align} with
\begin{align*}\hspace{-15mm} \begin{split}
\cF_m(\bsigma,\bu,A,H) &= \left\| \rho_{ice} H \frac{\partial \bu}{\partial t} + \bF(\bu) - \div\ {\bsigma} \right\|_0^2, \\
\quad \cF_e(\bsigma,\bu,A,H) &= \left\| \frac{\partial H}{\partial t}+\div(\bu H) \right\|_0^2+
\left\|  \frac{\partial A}{\partial t}+\div(\bu A)\right\|_0^2, \\
\cF_c(\bsigma,\bu,A,H) &= 
\left\| \bsigma- \frac {P(A,H)} 2 \left( \frac {\dev \ \bepsilon(\bu)} {\Delta(\bu)} + \frac  {\tr\ \bepsilon(\bu)} {\Delta(\bu)} \bI -\bI \right)\right\|_0^2 \ .\\
\end{split}\end{align*}
The time discretization can be realised using a $\theta$-scheme and decoupling the advection equations from the rest of the system such that for each time step $n+1$, the linear functional 
\begin{align}\begin{split}
\hspace{-15mm}\cG^{n+1} (A^{n+1},H^{n+1};\bu^{n},H^n,A^n)
=&\left\| \frac{H^{n+1}-H^n}{t^\Delta}+\div(\bu^{n} H^{n+1}) \right\|_0^2\\&+
\left\|  \frac{A^{n+1}-A^n}{t^\Delta}+\div(\bu^{n} A^{n+1})\right\|_0^2
\end{split}\end{align}
is first minimized over all $(A^{n+1},H^{n+1}) \in \left(L^2(\Omega)\right)^2$,
and then the functional
\begin{align}\begin{split}
\hspace{-15mm}\cF^{n+1}(\bsigma^{n+1}&,\bu^{n+1};\bsigma^{n},\bu^{n},A^{n+1},H^{n+1})\\
&= \left\| \rho_{ice} H^{n+1} \frac{\bu^{n+1}-u^n}{t^\Delta} + \bF(\bu^{n+\theta}) - \div\ \bsigma^{n+\theta}\right\|_0^2 \\
&\ +\cF_c (\bsigma^{n+1},\bu^{n+1};A^{n+1},H^{n+1})
\end{split}\end{align}
with the time discretized variables $$\bu^{n+\theta} = \theta\bu^{n+1} +(1-\theta)\bu^{n} $$ and
$$\bsigma^{n+\theta} = \theta\bsigma^{n+1} +(1-\theta)\bsigma^{n}\ , $$
is minimized over all $(\bsigma^{n+1},\bu^{n+1}) \in \left(H_{\text{div}}(\Omega)\right)^2 \times \left(H^1_{\Gamma_D}(\Omega)\right)^2$. For the spacial discretization, a conforming subspace $\bW_h$ of $\left(H_{\text{div}}(\Omega)\right)^2 \times \left(H^1_{\Gamma_D}(\Omega)\right)^2\times \left(L^2(\Omega)\right)^2$. Therefore, a triangulation $\mathcal T_h$ of the domain $\Omega$ is considered. In this work, we choose $\bW_h = (RT_1^2(\mathcal T_h) \times \mathcal P_2^2(\mathcal T_h)) \times \mathcal P_1^2(\mathcal T_h)) $ in order to have appropriate convergence properties.

%For the time-discretization, we introduce a discretization $\mathcal T_t = \{ [t_0,t_1],[t_1,t_2],...,[t_{n-1},t_n]\} $ with $t_0 = 0$ and $t_n=T$ of the inverval $[0,T]$. Note that in \cite{MehRic:17}, implicit Euler time stepping method the momentum equation and the balance law are decoupled in order to cope with the complex system. At this place, further investigations are needed. In this work however, we consider the coupled system and test the discretization for several $\theta$-scheme. At each time step $t_i$, the finite difference discretization of the time derivatives therefore reads
%\begin{align}\begin{split}
%%
%\forall \tau \in H_{div}(\Omega)
%%
%\forall \bv \in H^1_0(\Omega)
%%
%\forall p \in L^2(\Omega)
%%
%\forall q \in L^2(\Omega)
%\end{split}\end{align}
For the minimization of the nonlinear Functional $\mathcal \cF^{n+1}$ in each time step, the Least-Squares Functional is linearized
around a given approximation $(\bsigma^k,\bu^k,A^k,H^k)$ and the minimization is then carried out iteratively solving a sequence of linearized least squares problems. Additionaly the Least-Squares Functional is minimized subejct to the linear inequality constraints $A \in [0,1]$ and $H\geq 0$, that leads to a constraint optimization problem that we solved with an active set strategy.
% This constraint optimization problem may be solved, for example, by 
%Ellipticity future work.  \cite{BofBreFor:13}
Since the variables $A$ and $H$ are now decoupled from $\bu$ and $\bsigma$, we can define the stress-strain relation ship by
 \begin{align}\begin{split} 
\mathcal C(\bu;A,H):=\bsigma(\bu;A,H)  = \frac {P(A,H)} 2 
	\left(
		\frac { {\dev \ \bepsilon(\bu)}+  2 {\tr\ \bepsilon(\bu)}  \bI} {\Delta(\bu)}  -\bI 
	 \right) 
\end{split}\end{align}
The Gateaux derivative of $\mathcal C(\bu;A,H)$ in direction $\bv$ is denoted by $\mathcal C(\bu;A,H)[\bv]$ and given by
\begin{align}\begin{split}
\mathcal{J_C}(\bu;A,H)[\bv] = \frac {P(A,H)} 2  \left(
		\frac { {\dev \ \bepsilon(\bv)}+  2 {\tr\ \bepsilon(\bv)}  \bI} {\Delta(\bu)}
		+   \mathcal{J}_{\Delta^{-1}}(\bu)[\bv] \left(
		{\dev \ \bepsilon(\bu)}+  2 {\tr\ \bepsilon(\bu)}  \bI
		\right)
	 \right) 
\end{split}\end{align}
with $
  \mathcal{J}_{\Delta^{-1}}(\bu)[\bv] = -\Delta(\bu)^{-3} 
\left( \dev\ \bepsilon(\bu):\dev\ \bepsilon(\bv)+{4}\tr(\bepsilon(\bu))\tr(\bepsilon(\bv)) \right)
$. %For the computation of the first variation of the minimization of $\mathcal \cF^{n+1}$, the Gateaux derivative of $\bF(\bu)$, i.e. the Gateaux derivative of $\btau_o$ is needed:
%\begin{align}\begin{split}
%\mathcal{J_{\btau_o}}(\bu)[\bv] = \frac {P(A,H)} 2  \left(
%		\frac { {\dev \ \bepsilon(\bv)}+  2 {\tr\ \bepsilon(\bv)}  \bI} {\Delta(\bu)}
%		+   \mathcal{J}_{\Delta^{-1}}(\bu)[\bv] \left(
%		{\dev \ \bepsilon(\bu)}+  2 {\tr\ \bepsilon(\bu)}  \bI
%		\right)
%	 \right) 
%\end{split}\end{align}
%with $
%  \mathcal{J}_{\Delta^{-1}}(\bu)[\bv] = -\Delta(\bu)^{-3} 
%\left( \dev\ \bepsilon(\bu):\dev\ \bepsilon(\bv)+{4}\tr(\bepsilon(\bu))\tr(\bepsilon(\bv)) \right)
%$.

The first variation of the minimization of $\mathcal \cF^{n+1}$ is then given by
\begin{align*}\hspace{-15mm}\begin{split}
\mathcal B(\bsigma^{n+1},&\bu^{n+1};\btau,\bv;\bsigma^{n},\bu^{n},A^{n+1},H^{n+1})  = 
\left. \frac{\partial  \cF^{n+1}(\bsigma^{n+1}+\tau \btau,\bu^{n+1}+\tau\bv;\bsigma^{n},\bu^{n},A^n,H^n)}{\partial \tau}
\right|_{\tau = 0} % \\&
%\left. 
%\frac{\partial}{\partial \tau}
%\left( \bsigma^{n+1}+\tau \btau - \mathcal C (\bu^{n+1}+\tau\bv;A^{n+1},H^{n+1}), \bsigma^{n+1}+\tau \btau 
%- \mathcal C (\bu^{n+1}+\tau\bv;A^{n+1},H^{n+1}) \right)
%\right|_{\tau = 0} 
%\\&
%+
%\left. 
%\frac{\partial}{\partial \tau}
%\left(
%\rho_{ice} H^{n+1} \frac{\bu^{n+1}+\tau \bv-u^n}{t^\Delta} 
%+ \bF(\theta\bu^{n+1}+\tau \bv +(1-\theta)\bu^{n}) 
%- \div( \theta\bsigma^{n+1} +\tau \btau+(1-\theta)\bsigma^{n} ), 
%\rho_{ice} H^{n+1} \frac{\bu^{n+1}+\tau \bv-u^n}{t^\Delta} 
%+ \bF(\theta\bu^{n+1}+\tau \bv +(1-\theta)\bu^{n}) 
%- \div( \theta\bsigma^{n+1} +\tau \btau+(1-\theta)\bsigma^{n} )
%\right)
%\right|_{\tau = 0}
 \\=&
\left. 
\frac{\partial}{\partial \tau}
\left( \bsigma^{n+1}+\tau \btau - \mathcal C (\bu^{n+1}+\tau\bv;A^{n+1},H^{n+1}), \bsigma^{n+1}+\tau \btau 
- \mathcal C (\bu^{n+1}+\tau\bv;A^{n+1},H^{n+1}) \right)
\right|_{\tau = 0} \\&
+
2
\left(
\rho_{ice} H^{n+1} \frac{\bu^{n+1}-u^n}{t^\Delta} 
- \div( \theta\bsigma^{n+1}+(1-\theta)\bsigma^{n} ), 
\rho_{ice} H^{n+1} \frac{ \theta \bv}{t^\Delta} 
- \div(\theta \btau )
\right)
\\&
+
\left. 
\frac{\partial}{\partial \tau}
\left(
\bF(\theta\bu^{n+1}+\theta \tau \bv +(1-\theta)\bu^{n}), 
\bF(\theta\bu^{n+1}+\theta \tau \bv +(1-\theta)\bu^{n}) 
\right)
\right|_{\tau = 0}
 \\+&
2\left. 
\frac{\partial}{\partial \tau}
\left(
\rho_{ice} H^{n+1} \frac{\bu^{n+1}+\tau \bv-u^n}{t^\Delta} 
- \div( \theta\bsigma^{n+1} +\theta\tau \btau+(1-\theta)\bsigma^{n} ), 
\bF(\theta\bu^{n+1}+\theta\tau \bv +(1-\theta)\bu^{n}) 
\right)
\right|_{\tau = 0}
 \\=&
2 \left( 
\bsigma^{n+1} - \mathcal C (\bu^{n+1};A^{n+1},H^{n+1}), 
\btau 
- J_{\mathcal C} (\bu^{n+1};A^{n+1},H^{n+1})[\bv] \right) \\
&+2 \left(\rho_{ice} H^{n+1} \frac{\bu^{n+1}-u^n}{t^\Delta} + \bF(\bu^{n+\theta}) - \div\ \bsigma^{n+\theta} ,
\rho_{ice} H^{n+1} \frac{\bv}{t^\Delta} + J_{\bF}(\bu^{n+\theta})[\theta \bv] - \theta \div\ \btau \right) \ .
\end{split}\end{align*}
Setting this first variation to zero leads to a necessary condition such that the Gau\ss -Newton Methods in each time step consits in setting iterativly  
$
(\bsigma^{n+1,k+1}, \bu^{n+1,k+1})=(\bsigma^{n+1,k}, \bu^{n+1,k})+(\delta \bsigma, \delta \bu)
$
where $(\delta \bsigma, \delta \bu) \in (RT_0^2(\mathcal T_h) \times \mathcal P_1^2(\mathcal T_h))$ is the solution of 
\begin{align*}\hspace{-15mm} \begin{split}
\mathcal B(&\bsigma^{n+1,k},\bu^{n+1,k};\btau,\bv;\bsigma^{n},\bu^{n},A^{n+1},H^{n+1}) =  \left( 
\delta \bsigma 
- J_{\mathcal C} (\bu^{n+1,k};A^{n+1},H^{n+1})[\delta \bu]  , 
\btau 
- J_{\mathcal C} (\bu^{n+1,k};A^{n+1},H^{n+1})[\bv] \right) \\&
+ \left(
\rho_{ice} H^{n+1} \frac{\delta \bu}{t^\Delta} + J_{\bF}(\theta\bu^{n+1,k}+(1-\theta)\bu^n)[\theta \delta \bu] - \theta \div\ \delta\bsigma ,
\rho_{ice} H^{n+1} \frac{\bv}{t^\Delta} + J_{\bF}(\theta\bu^{n+1,k}+(1-\theta)\bu^n)[\theta \bv] - \theta \div\ \btau
 \right)\\&
\text{for all $(\btau,\bsigma) \in (RT_0^2(\mathcal T_h) \times \mathcal P_1^2(\mathcal T_h))$.}\end{split}\end{align*}

\section{Test Case}
\begin{figure}[t]
\includegraphics[height=5cm]{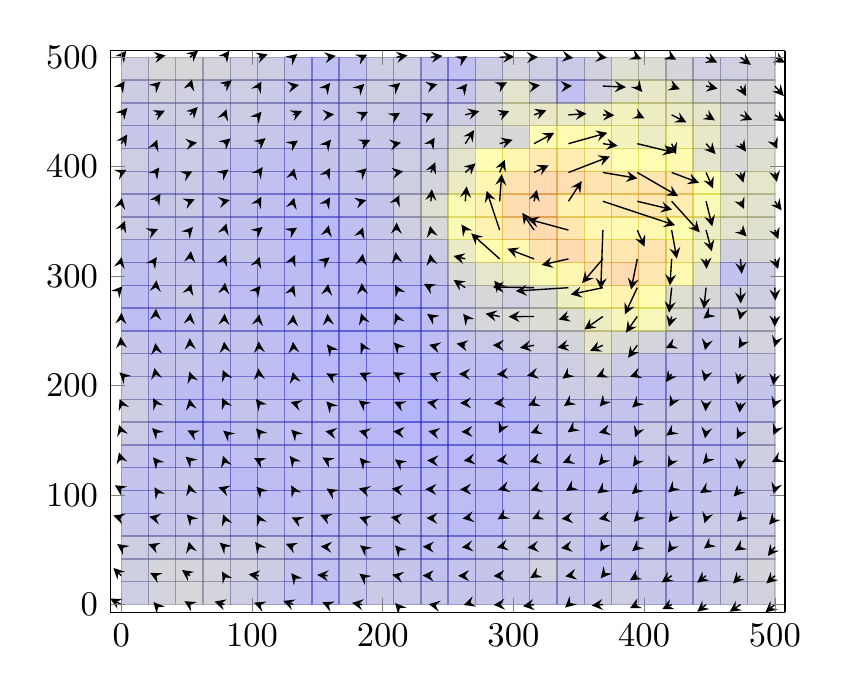}
\hfil
\includegraphics[height=5cm]{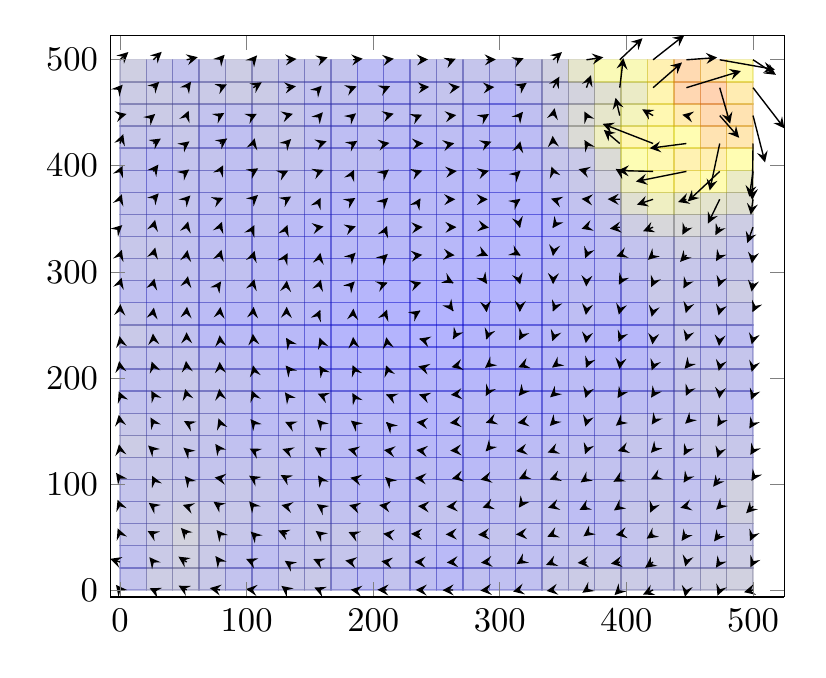}\\
\includegraphics[height=5cm]{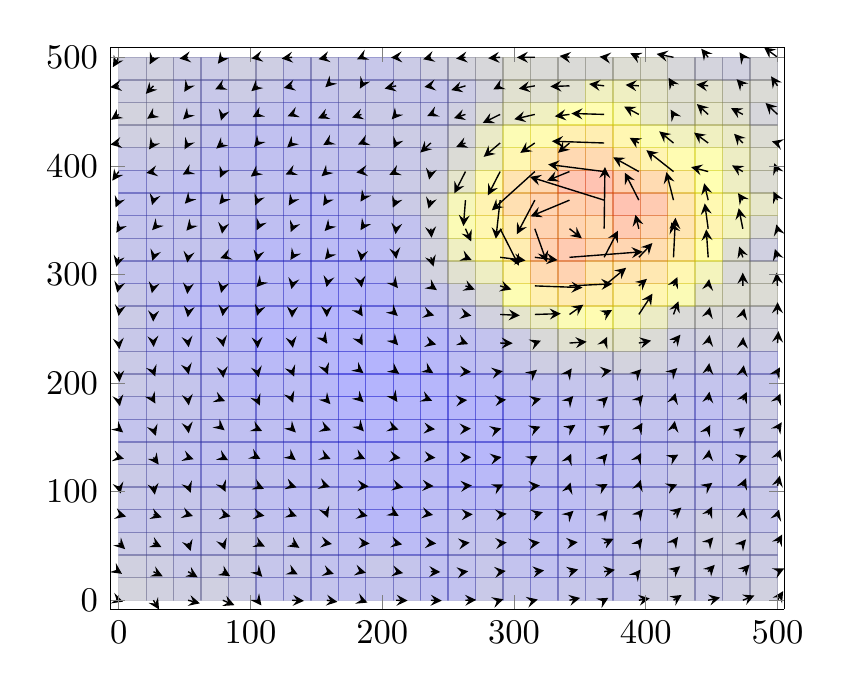}
\hfil
\includegraphics[height=50mm]{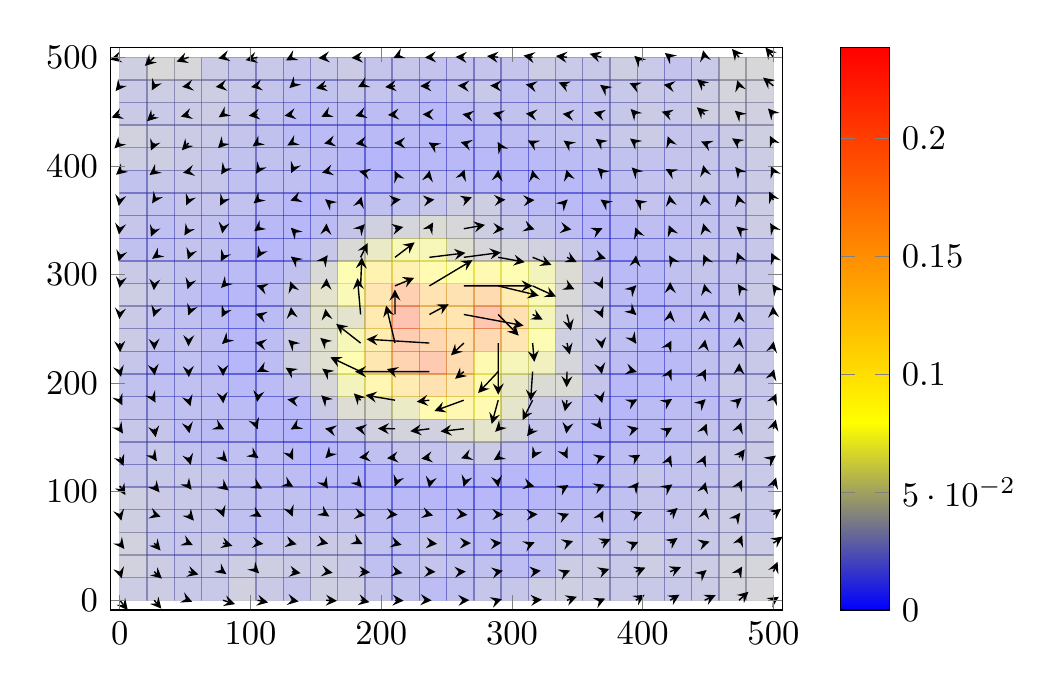}\\
\caption{Sea-ice velocity at $t=2,4,6,8$.}
\label{fig:5}
\end{figure}

In order to investigate the approximation properties of the Least-Squares method, we consider the same test case as in \cite{MehRic:17}, involving a quadratic domain (see also \cite{Hunke:01}) and simulating the sea ice dynamics for $T=8$ days. Since the Least-Squares Method approximates all the residuals of the partial differential equation simultaneously, we scale the domain to the unit square $\Omega = [0,1]^2$. Since the wind field is a cyclone from the midpoint of the computational domain to the edge followed by an anticyclone diagonally passing from the edge to the midpoint, we define the time $t^m = t-4$ measured in days with respect to the time when the wind forcing alternates from cyclonic to anticyclonic. Further, let $\tilde\bx(t) = \bx-\bx^m(t)$ denote the position with respect to the center of the cyclone $\bx^m(t)=x^m(t)(\be_1+\be_2)$ with $x^m(t)= 0.1(9-|t^m|) $.  Then, the prescribed wind field is given by
\begin{align}
\bv_a = 10 {v_a^m} \left( 1-\frac{2}{e^{t^m}e^{8-|t^m|}+1}\right) e^{-\frac{\|\tilde\bx(t)\|_2}{10}} 
\bR\left(\frac{17}{40}\pi+\frac{t^m}{40|t^m|}\pi\right) \tilde \bx  & \\
\quad \text{with } \bR(\vartheta) =  {\begin{pmatrix}\cos \vartheta &-\sin \vartheta \\ \sin \vartheta &\cos \vartheta \end{pmatrix}} \ ,
\end{align}
	and a maximal wind velocity $v_a^m$, while the circular steady ocean current is
	\begin{align}
	\bv_o = v_o^m \begin{pmatrix} 2y-1 \\ 1-2x \end{pmatrix}
	\end{align}
	with a maximal ocean velocity $v_o^m$. Finally, the initial conditions are given by zero velocity, constant ice concentration $A=1$ and $H^0(x,y)=0.3+0.005(\sin(250x)+\sin(250y))$. All simulations are executed with Fenics, using the inherent Newton solver. The velocity results at $t=2,4,6,8$ days are shown in the figure \ref{fig:5}. Further intervestigations are needed, in particular regarding the ellipticity of the Least-Squares Functional, the possibility of considering domain with curved boundaries (see \cite{BerMueSta:14b}) and the relation to others standard or mixed methods (as in \cite{Brandts}).

\begin{figure}
	\begin{tabular}{@{}lll@{}}
	\hline
	Parameter & Value  \\
	\hline
	maximal ocean velocity $v_o^m$ & $0.01$ ms$^-1$ \\
	maximal ocean velocity $v_a^m$ & $15$ ms$^-1$  \\
	sea ice density $\rho_{ice}$ & $900$ kg\ m$^-3$  \\
	air density $\rho_a$ & $1.3$ kg\ m$^-3$  \\
	water density $\rho_o$ & $1026$ kg\ m$^-3$  \\
	air drag coefficient  $C_a$ & $1.2 \cdot 10^{-3}$  \\
	water  $C_o$&  $5.5 \cdot 10^{-3}$  \\
	coriolis parameter $f_c$	& $1.46 \cdot 10^{-4}$ s$^-1$  \\
	ice strength parameter  $P^\star$&  $27.5 \cdot 10^{3}$ Nm$^-2$  \\
	ice concentration parameter $C$ & 20\\
	\hline
	\end{tabular}
	\caption{Parameter used in the simulation}
	\label{tab:4}
	\end{figure}

%As in \cite{} and \cite{}

\end{document}